\documentclass{amsart}
\DeclareMathSymbol{\twoheadrightarrow}
{\mathrel}{AMSa}{"10}

\def\Q{{\mathbf Q}}
\def\Z{{\mathbf Z}}
\def\C{{\mathbf C}}
\def\R{{\mathbf R}}
\def\F{{\mathbf F}}
\def\R{{\mathfrak R}}
\def\bchi{{\mathbf \chi}}
\def\bphi{{\mathbf \phi}}
\def\CC{{\mathfrak C}}

\def\SS{{\mathbf S}}
\def\A{{\mathbf A}}
\def\Sn{{\mathbf S}_n}
\def\An{{\mathbf A}_n}

\def\Gal{\mathrm{Gal}}
\def\Perm{\mathrm{Perm}}
\def\Alt{\mathrm{Alt}}

\def\End{\mathrm{End}}
\def\Aut{\mathrm{Aut}}

\def\Mat{\mathrm{Mat}}
\def\cl{\mathrm{cl}}

\def\I{\mathrm{Id}}

\def\fchar{\mathrm{char}}

\def\GL{\mathrm{GL}}
\def\SL{\mathrm{SL}}
\def\PSL{\mathrm{PSL}}

\def\PGL{\mathrm{PGL}}

\def\dim{\mathrm{dim}}

\def\P{{\mathbf P}}

\newtheorem{thm}{Theorem}[section]
\newtheorem{lem}[thm]{Lemma}
\newtheorem{cor}[thm]{Corollary}

\theoremstyle{definition}
\newtheorem{defn}[thm]{Definition}
\newtheorem{ex}[thm]{Example}

\newtheorem{rem}[thm]{Remark}
\newtheorem{rems}[thm]{Remarks}
\title[Cyclic covers, jacobians and endomorphisms]
{Cyclic covers of the projective line, their jacobians and endomorphisms}
\author[Yuri G. Zarhin]{Yuri G. Zarhin}
\address{Department of Mathematics, Pennsylvania State University,
University Park, PA 16802, USA}
\email{zarhin\char`\@math.psu.edu}
\thanks{Partially supported by the NSF}
\begin{document}
\maketitle
\section{Introduction}
As usual, $\Z,\Q,\C$ denote the ring of integers, the field of rational numbers and the field of complex numbers respectively.
Let $p$ be an odd prime. Recall that $p$ is called
a {\sl Fermat prime} if $p=2^{2^r}+1$ for some
positive integer $r$; e.g., $3,5,17,257$ are Fermat
prime.

Let us fix a primitive $p$th root of unity
$$\zeta_p \in \C.$$
Let $\Q(\zeta_p)$ be the $p$th cyclotomic field.
It is well-known that $\Q(\zeta_p)$ is a CM-field.
If $p$ is a Fermat prime then the only CM-subfield
of $\Q(\zeta_p)$ is $\Q(\zeta_p)$ itself, since
the Galois group of $\Q(\zeta_p)/\Q$ is a cyclic $2$-group, whose only element of order $2$ acts
as the complex conjugation. All other subfields of
$\Q(\zeta_p)$ are totally real.

Let $f(x) \in \C[x]$ be a polynomial of degree $n\ge 4$ without multiple roots. Let $C_{f,p}$ be a smooth projective model of the smooth affine curve
$$y^p=f(x).$$
It is well-known
(\cite{Koo}, pp. 401--402, \cite{Towse}, Prop. 1 on p. 3359,
 \cite{Poonen}, p. 148)
that the genus $g(C_{f,p})$  of $C_{f,p}$ is $(p-1)(n-1)/2$ if $p$ does not divide $n$ and $(p-1)(n-2)/2$ if it does.
The map
 $$(x,y) \mapsto (x, \zeta_p y)$$
gives rise to a non-trivial birational autmorphism
$$\delta_p: C_{f,p} \to C_{f,p}$$
of period $p$.

Let $J^{(f,p)}= J(C_{f,p})$  be the jacobian of $C_{f,p}$;
it is an abelian variety, whose dimension equals $g(C_{f,p})$.
We write $\End(J^{(f,p)})$ for the ring of endomorphisms of $J^{(f,p)}$. By functoriality,
$\delta_p$ induces an automorphism of $J^{(f,p)}$
which we still denote by $\delta_p$; it is known (\cite{Poonen}, p. 149, \cite{SPoonen}, p. 448)) that
$$\delta^{p-1}+\cdots +\delta_p+1=0$$
in $\End(J^{(f,p)})$. This gives us an embeddding
$$\Z[\zeta_p] \cong \Z[\delta_p] \subset \End(J^{(f,p)})$$
(\cite{Poonen}, p. 149, \cite{SPoonen}, p. 448)).

Our main result is the following statement.

\begin{thm}
\label{endo}
Let $K$ be a subfield of $\C$ such that all the coefficients of $f(x)$
lie in $K$. Assume also that $f(x)$ is an irreducible polynomial in $K[x]$
of degree $n \ge 5$
 and its Galois group over $K$ is either the symmetric group $\Sn$ or
the alternating group $\An$.
Then
$\Z[\delta_p]$ is a maximal commutative subring in
$\End(J^{(f,p)})$.

If $p$ is a a Fermat prime (e.g., $p=3,5,17,257$) then
$$\End(J^{(f,p)})=\Z[\delta_p] \cong\Z[\zeta_p].$$
\end{thm}

 When $p=3$ one may obtain an additional information about
Hodge classes on self-products of the corresponding trigonal jacobian.

\begin{thm}
\label{trigonal}
Let $K$ be a subfield of $\C$ such that all the coefficients of $f(x)$
lie in $K$. Assume also that $f(x)$ is an irreducible polynomial in $K[x]$
 of degree $n \ge 5$
 and its Galois group over $K$ is either the symmetric group $\Sn$ or
the alternating group $\An$. If $3$ does not divide $n-1$ then:
\begin{enumerate}
\item[(i)]
Every Hodge class
on each self-product of $J^{(f,3)}$ could be presented as a linear
combination of products of divisor classes.
In particular, the Hodge conjecture is true for each self-product of
 $J^{(f,3)}$.

\item[(ii)]
If $K$ is a number field containing $\sqrt{-3}$ then every Tate class on
each self-product of $J^{(f,3)}$ could be presented as a linear combination
of products of divisor classes. In particular, the Tate conjecture is true
for each self-product of $J^{(f,3)}$.
\end{enumerate}
\end{thm}

\begin{ex}
The polynomial $x^n-x-1 \in \Q[x]$ has Galois group $\Sn$ over $\Q$
(\cite{Serre}, p. 42). Therefore the ring of endomorphism (over $\C$)
of the jacobian $J(C^{(n,3)})$ of the  curve $C^{(n,3)}:y^3=x^n-x-1$ is
 $\Z[\zeta_3]$ if  $n\ge 5$.

If  $n=3k-1$ for some integer $k \ge 2$ then
all Hodge classes on each self-products of $J(C^{(n,3)})$  could be
presented as linear combinations of products of divisor classes.
In particular, the Hodge conjecture is true for all these self-products.
Notice that $J(C^{(n,3)})$ is an abelian variety defined over $\Q$ of
dimension $n-1=3k-2$.
\end{ex}

\begin{rems}
\label{ribcom}
\begin{enumerate}
\item[(i)]
If $f(x) \in K[x]$ then the curve
$C_{f,p}$ and its jacobian $J^{(f,p)}$ are defined over $K$. Let $K_a\subset \C$ be the algebraic closure of $K$. Clearly, all endomorphisms of $J^{(f,p)}$ are defined over $K_a$.
 This implies that in order to prove Theorem \ref{endo}, it suffices to
check that $\Z[\delta_p]$ is a maximal commutative subring in the ring of $K_a$-endomorphisms of $J^{(f,p)}$ or equivalently, that
$\Q[\delta_p]$ is a maximal commutative $\Q$-subalgebra in the algebra of $K_a$-endomorphisms of $J^{(f,p)}$.

\item[(ii)]
Assume that $p=3$ and $\Z[\delta_3]=\End(J^{(f,3)})$.
 The endomorphism algebra
$\End^0(J^{(f,p)})=\End(J^{(f,p)})\otimes\Q$
is the imaginary quadratic field
$\Q(\delta_3)\cong\Q(\sqrt{-3})$.

 There are exactly two embeddings
$$\sigma, \bar{\sigma}:\Q(\delta_3)\hookrightarrow K_a\subset\C$$
and they are complex-conjugate. We have
$$\Q(\delta_3)\otimes_{\Q}\C=\C \oplus \C.$$
By functoriality, $\Q(\delta_3)$
acts on the $\C$-vector space
$H^{1,0}(J^{(f,3)})=\Omega^1(J^{(f,3)})$ of differentials of the fist kind.
This action gives rise to a splitting of the
$\Q(\delta_3)\otimes_{\Q}\C=\C \oplus \C$-module
$$H^{1,0}(J^{(f,3)})=H^{1,0}_{\sigma}\oplus H^{1,0}_{\bar{\sigma}}.$$
The dimensions $n_{\sigma}:=\dim_{\C}(H^{1,0}_{\sigma})$
and $n_{\bar{\sigma}}:=\dim_{\C}(H^{1,0}_{\bar{\sigma}})$
are called the {\sl multiplicities} of $\sigma$ and $\bar{\sigma}$
respectively. Clearly, $n_{\sigma}$  (resp. $n_{\bar{\sigma}}$) coincides
with the multiplicity of the eigenvalue $\sigma(\delta_3)$
(resp. $\bar{\sigma}(\delta_3)$) of the induced $\C$-linear operator
$$\delta_3^*:\Omega^1(J^{(f,3)})\to\Omega^1(J^{(f,3)}).$$
By a theorem of Ribet (\cite{Ribet3}, Th. 3 on p. 526),
 if the multiplicities $n_{\sigma}$
and $n_{\bar{\sigma}}$ are {\sl relatively prime} and
$\End^0(J^{(f,3)})=\Q(\delta_3)$ then every Hodge class on each
self-product of $J^{(f,3)}$ could be presented as a linear combination
of products of divisor classes.
Therefore, the assertion (i) of Theorem \ref{trigonal} would follow
 from Theorem \ref{endo} (with $p=3$) if we know that the multiplicities
$n_{\sigma}$ and $n_{\bar{\sigma}}$ are relatively
prime while $3$ does not divide $n-1$.
\item[(iii)]
One may easily check that $n_{\sigma}$ (resp. $n_{\bar{\sigma}}$) coincides
 with the multiplicity of the eigenvalue $\sigma(\delta_3)$
(resp. $\bar{\sigma}(\delta_3)$) of the induced $\C$-linear operator
$$\delta_3^*:\Omega^1(C_{(f,3)})\to\Omega^1(C_{(f,3)}).$$
\end{enumerate}
\end{rems}

\section{Permutation groups and permutation modules}
\label{permute}

Let $B$ be a finite set consisting of $n \ge 5$ elements.
We write $\Perm(B)$ for the group of permutations of $B$.
A choice of ordering on $B$ gives rise to an isomorphism
$$\Perm(S) \cong \Sn.$$
We write $\Alt(B)$ for the only subgroup in $\Perm(B)$ of index $2$. Clearly, $\Alt(B)$ is normal and isomorphic to the alternating group $\An$. It is well-known that $\Alt(B)$ is a simple non-abelian group of order $n!/2$.
Let $G$ be a  subgroup of $\Perm(B)$.

Let $\F$ be a field. We write $\F^B$
for the $n$-dimensional $\F$-vector space of maps $h:B \to F$.
The space $\F^B$ is  provided with a natural action of $\Perm(B)$ defined
as follows. Each $s \in \Perm(B)$ sends a map
 $h:B\to \F_2$ into  $sh:b \mapsto h(s^{-1}(b))$. The permutation module $\F^B$
contains the $\Perm(B)$-stable hyperplane
$$(\F^B)^0=
\{h:B\to\F\mid\sum_{b\in B}h(b)=0\}$$
and the $\Perm(B)$-invariant line $\F \cdot 1_B$ where $1_B$ is the constant function $1$. The quotient $\F^B/(\F^B)^0$ is a trivial $1$-dimensional $\Perm(B)$-module.

Clearly, $(\F^B)^0$ contains $\F \cdot 1_B$ if and only if $\fchar(\F)$ divides $n$.
 If this is {\sl not} the case then there is a $\Perm(B)$-invariant splitting
$$\F^B=(\F^B)^0 \oplus \F \cdot 1_B.$$

Clearly, $\F^B$ and $(\F^B)^0$  carry natural structures of $G$-modules.
Their (Brauer) characters depend only on characteristic of $F$.

Let us consider the case of $F=\Q$. Then
 the character of $\Q^B$ sends each $g \in G$ into the number of fixed points of $g$ (\cite{SerreRep}, ex. 2.2,p.12); it takes on values in $\Z$ and called the
 {\sl permutation character} of $B$. Let us denote by
$\bchi=\bchi_B:G \to \Q$
 the character of $(\Q^B)^0$.

It is known that the $\Q[G]$-module $(\Q^B)^0$ is absolutely simple if and only if $G$ acts doubly transitively on $B$ (\cite{SerreRep}, ex. 2.6, p. 17).
 Clearly, $1+\bchi$ is the permutation character. In particular, $\bchi$ also takes on values in $\Z$.

Now, let us consider the case of $\F=\F_p$.

If  $p\mid n$  then let us define the $\Perm(B)$-module
$$(\F_p^B)^{00}:=(\F_p^B)^0/(\F_p \cdot 1_B).$$
If $p$ does not divide $n$ then let us put
$$(\F_p^B)^{00} := (\F_p^B)^0.$$

\begin{rem}
Clearly,
$\dim_{\F_p}((\F_p^B)^{00})=n-1$ if $n$ is not divisible by $p$
 and $\dim_{\F_p}((\F_p^B)^{00})=n-2$ if $p\mid n$.
In both cases $(\F_p^B)^{00}$ is a faithful $G$-module.
One may easily check that if the $\F_p[G]$-module $(\F_p^B)^{00}$
is absolutely simple then the $\Q[G]$-module $(\Q^B)^0$ is
also absolutely simple and therefore $G$ acts doubly transitively on $B$.
\end{rem}

 Let $G^{(p)}$ be the set of $p$-regular elements of $G$. Clearly,
the Brauer character of the $G$-module $\F_{p}^B$
coincides with the restriction of $1+\bchi_B$ to $G^{(p)}$.
This implies easily that the Brauer character of the $G$-module $(\F_{p}^B)^0$
coincides with the restriction of  of $\bchi_B$ to $G^{(p)}$.

\begin{rem}
\label{Bcharacter}
 Let us denote by
$\bphi_B=\bphi$
 the Brauer character of the $G$-module $(\F_p^B)^{00}$.
 One may easily check that $\bphi_B$ coincides with the restriction of
$\bchi_B$ to $G^{(p)}$ if $p$ does not divide $n$  and
 with the restriction of $\bphi_B-1$ to $G^{(p)}$ if  $p\mid n$.
In both cases $\bphi_B$ takes on values in $\Z$.
\end{rem}

\begin{ex}
\label{A5p5}
Suppose $n=p=5$ and $G=\Alt(B)\cong \A_5$.
Then in the notations of \cite{Atlas}, p. 2 and \cite{AtlasB}, p. 2
$\bchi_B=1+\chi_4$ and the restriction of $\bphi_B-1$=$\chi_4-1$ to $G^{(5)}$ coincides with absolutely irreducible Brauer character $\varphi_2$. This implies that the $\Alt(B)$-module $(\F_p^B)^{00}$ is absolutely simple.
\end{ex}

The following elementary assertion is based on Lemma 7.1 on p. 52 of \cite{Passman} and Th. 9.2 on p. 145 of \cite{Isaacs}. (The case of $p=2$ is Lemma 5.1 of \cite{Zarhin}).

\begin{lem}
\label{semisimple}
 Assume that $G$ acts doubly transitively on $B$. If $p$ does not divide $n$ then $\End_{G}((\F_p^{B})^0)=\F_p$. In particular, if the $G$-module $(\F_p^{B})^0$ is semisimple then it is absolutely simple.
\end{lem}

\begin{proof}
It suffices to check that
$\dim_{\F_p}(\End_{G}((\F_p^{B})^0)) \le 1$. Itn order to do that, recall that the double transitivity implies that
$\dim_{\F_p}(\End_{G}((\F_p^{B}))) =2$
(Lemma 7.1 on p. 52 of \cite{Passman}). Now the desired inequality
follows easily from the existence of the $G$-invariant splitting
$$\F_p^B=(\F_p^B)^0 \oplus \F_p \cdot 1_B.$$
\end{proof}

\begin{rem}
\label{oddeven}
Assume that $n=\#(B)$ is divisible by $p$.
Let us choose $b \in B$ and let $G':=G_b$ be the stabilizer of $b$ in $G$
 and $B'=B\setminus\{b\}$.
Then $n'=\#(B')=n-1$ is not divisible by $p$ and there is a canonical isomorphism
 of $G'$-modules
$$(\F_p^{B'})^{00}   \cong (\F_p^B)^{00}$$
defined as follows.
First, there is a natural $G'$-equivariant embedding $\F_p^{B'} \subset \F_p^B$
which could be
obtained by extending each $h: B' \to \F_p$ to $B$ by letting $h(b)=0$.
Second, this embedding identifies $(\F_p^{B'})^0$ with a hyperplane of
$(\F_p^{B})^0$ which does not contain $1_B$. Now the composition
$$(\F_p^{B'})^{00}  =(\F_p^{B'})^0 \subset (\F_p^{B})^0 \to
(\F_p^{B})^0/(\F_p \cdot 1_B)=(\F_p^B)^{00}$$
gives us the desired isomorphism.
This implies that if the $G_b$-module $(\F_p^{B'})^{00}$ is absolutely simple then the $G$-module $(\F_p^{B})^{00}$ is also absolutely simple.

For example, if $G=\Perm(B)$ (resp. $\Alt(B)$) then
$G_b=\Perm(B')$ (resp. $\Alt(B')$) and therefore
the $\Perm(B')$-modules (resp. $\Alt(B')$-modules
$(\F_p^{B})^{00}$ and $(\F_p^{B'})^{00}$ are isomorphic. We use this observation in order to prove the following statement.
\end{rem}

The following assertion goes back to Dickson.

\begin{lem}
\label{Anp}
Assume that $G=\Perm(B)$ or $\Alt(B)$.
Then the $G$-module $(\F_p^{B})^{00}$ is absolutely simple.
\end{lem}

\begin{proof}
In light of Example \ref{A5p5}, we may assume that $(n,p)\ne (5,5)$.
In light of Remark \ref{oddeven} we may assume that  $p$ does not divide $n$ and therefore $$(\F_p^{B})^{00}=(\F_p^{B})^0.$$
The natural representation of $\Perm(B)=\Sn$ in $(\F_p^{B})^0$ is irreducible
(\cite{Glasgow}, Th. 5.2 on p. 133).

Since $\Alt(B)$ is normal in $\Perm(B)$, the $\Alt(B)$-module $(\F_p^{B})^0$
is semisimple, thanks to Clifford's theorem
(\cite{CR}, \S 49, Th. 49.2).
 Since $n \ge 5$, the action of $\Alt(B)\cong \An$ on $B \cong \{1, \cdots, n\}$ is
doubly transitive. Applying Lemma \ref{semisimple}, we conclude that  the $\Alt(B)$-module $(\F_p^{\R_f})^0$ is absolutely simple. (See also \cite{Wagner}.)
\end{proof}

\section{Cyclic covers and jacobians}
\label{Prelim}

 Throughout this paper we fix a prime $p$ and assume
that $K$ is a field of characteristic
zero.
We fix
its algebraic closure $K_a$ and write $\Gal(K)$ for the absolute
Galois group $\Aut(K_a/K)$. We also fix in $K_a$ a primitive $p$th
root of unity $\zeta$.

Let $f(x) \in K[x]$ be a separable polynomial of degree $n \ge 4$.
We write $\R_f$ for the set of its roots and denote by
$L=L_f=K(\R_f)\subset K_a$ the corresponding splitting field. As
usual, the Galois group $\Gal(L/K)$ is called the Galois group of
$f$ and denoted by $\Gal(f)$. Clearly, $\Gal(f)$ permutes elements
of $\R_f$ and the natural map of $\Gal(f)$ into the group
$\Perm(\R_f)$ of all permutations of $\R_f$ is an embedding. We
will identify $\Gal(f)$ with its image and consider it as a
permutation group of $\R_f$. Clearly, $\Gal(f)$ is transitive if
and only if $f$ is irreducible in $K[x]$. Therefore the
$\Gal(f)$-module $(\F_p^{\R_f})^{00}$ is defined. The canonical
surjection $$\Gal(K) \twoheadrightarrow \Gal(f)$$ provides
$(\F_p^{\R_f})^{00}$ with canonical structure of the
$\Gal(K)$-module via the composition $$\Gal(K)\twoheadrightarrow
\Gal(f)\subset \Perm(\R_f) \subset \Aut((\F_p^{\R_f})^{00}).$$
Let us put $$V_{f,p}=(\F_p^{\R_f})^{00}.$$

 Let $C=C_{f,p}$ be the smooth projective model of the smooth
affine $K$-curve
            $$y^p=f(x).$$

So, $C$ is a smooth projective curve defined over $K$. The
rational function $x \in K(C)$ defined a finite cover  $\pi:C \to
\P^1$ of degree $p$. Let $B'\subset C(K_a)$ be the set of
ramification points.  Clearly, the restriction of $\pi$ to $B'$ is
an {\sl injective} map $\pi:B' \hookrightarrow \P^1(K_a)$, whose
image is the disjoint union of $\infty$ and  $\R_f$ if $p$ does
{\sl not} divide $\deg(f)$ and just $\R_f$ if it does.
We write
$$B=\pi^{-1}(\R_f)=\{(\alpha,0)\mid \alpha \in \R_f\} \subset B' \subset C(K_a).$$
Clearly, $\pi$ is ramified at each point of $B$ with ramification index $p$. We have $B'=B$ if and only if $n$
is  divisible by $p$. If $n$ is not divisible by $p$ then $B'$ is the disjoint union of $B$ and a single point
$\infty':=\pi^{-1}(\infty)$; in addition, the ramification index of $\pi$ at $\pi^{-1}(\infty)$ is also $p$. If $p$ does divide $n$ then $\pi^{-1}(\infty)$ consists of $p$ unramified points denoted by $\infty_1, \dots , \infty_p$. This implies that the inverse
image $\pi^*(n(\infty))=n \pi^*(\infty)$ of the divisor $n(\infty)$ is always divisible
 by $p$ in the divisor group of $C$.
  Using Hurwitz's formula, one may
easily compute genus $g=g(C)=g(C_{p,f})$ of $C$
(\cite{Koo}, pp. 401--402, \cite{Towse}, Prop. 1 on p. 3359,
 \cite{Poonen}, p. 148).
Namely, $g$ is $(p-1)(n-1)/2$ if $p$ does {\sl not} divide
$p$ and $(p-1)(n-2)/2$ if it does. See \S 1 of \cite{Towse}
for an explicit description of a smooth complete model of $C$ (when $n>p$).

 Assume that $K$ contains $\zeta$. There is a non-trivial birational
 automorphism of $C$
 $$\delta_p:(x,y) \mapsto (x, \zeta y).$$
Clearly, $\delta_p^p$ is the identity map and
 the set of fixed points of $\delta_p$ coincides with $B'$.

Let $J^{(f,p)}=J(C)=J(C_{f,p})$ be the jacobian of $C$. It is a $g$-dimensional abelian variety defined over $K$ and one may view $\delta_p$ as an element of
 $$\Aut(C) \subset\Aut(J(C)) \subset \End(J(C))$$
such that
  $$\delta_p \ne \I, \quad \delta_p^p=\I$$
where $\I$ is the identity endomorphism of $J(C)$.
Here $\End(J(C))$ stands for the ring of all $K_a$-endomorphisms of $J(C)$. As usual, we write $\End^0(J(C))=\End^0(J^{(f,p)})$ for
the corresponding $\Q$-algebra $\End(J(C))\otimes \Q$.

\begin{lem}
\label{cycl}
 $\I+\delta_p+ \cdots + \delta_p^{p-1}=0$ in $\End(J(C))$.
Therefore the subring $\Z[\delta_p] \subset \End(J(C))$ is
isomorphic to the ring $\Z[\zeta_p]$ of integers in the $p$th
cyclotomic field $\Q(\zeta_p)$. The $\Q$-subalgebra

$$\Q[\delta_p]\subset\End^0(J(C))=\End^0(J^{(f,p)})$$ is
isomorphic to $\Q(\zeta_p)$.
\end{lem}

\begin{proof}
See \cite{Poonen}, p. 149, \cite{SPoonen}, p. 448.
\end{proof}

\begin{rems}
\label{nondiv}
\begin{enumerate}
\item[(i)]

Assume that $p$ is odd and $n=\deg(f)$ is divisible by $p$ say, $n=pm$ for some
positive integer $m$. Then $n \ge 5$.

 Let $\alpha \in K_a$ be a root of $f$ and $K_1=K(\alpha)$ be
the corresponding subfield of $K_a$. We have $$f(x)=(x-\alpha)
f_1(x)$$ with $f_1(x) \in K_1[x]$. Clearly, $f_1(x)$
 is a separable
polynomial over $K_1$ of odd degree $pm-1=n-1 \ge 4$. It is also
clear that the polynomials
$$h(x)=f_1(x+\alpha), h_1(x)=x^{n-1}h(1/x) \in K_1[x]$$
are separable of the same degree $pm-1=n-1\ge 4$.

The standard substitution $$x_1=1/(x-\alpha), y_1=y/(x-\alpha)^m$$
establishes a birational isomorphism between $C_{f,p}$ and a
superelliptic curve $$C_{h_1}: y_1^p=h_1(x_1)$$
(see \cite{Towse}, p. 3359). But
$\deg(h_1)=pm-1$ is {\sl not} divisible by $p$. Clearly, this isomorphism commutes with the actions of $\delta_p$. In particular, it induces an isomorphism of
$\Z[\delta_p]$-modules $J^{(f,p)}(K_a)$ and $J^{(h_1,p)}(K_a)$ which commutes with the action of $\Gal(K_1)$.

\item[(ii)]
Assume, in addition, that $f(x)$ is irreducible in $K[x]$ and
$\Gal(f)$ acts $s$-transitively on $\R_f$ for some positive integer
$s \ge 2$. Then the Galois group $\Gal(h_1)$ of $h_1$ over $K_1$
acts $s-1$-transitively on the set $\R_{h_1}$ of roots of $h_1$. In
particular, $h_1(x)$ is irreducible in $K_1[x]$.

It is also clear that if $\Gal(f)=\Sn$ or $\An$ then
$\Gal(h_1)=\SS_{n-1}$ or $\A_{n-1}$ respectively. \end{enumerate}
\end{rems}

Let us put $\eta=1-\delta_p$. Clearly, $\eta$ divides $p$ in
$\Z[\delta_p]\cong \Z[\zeta_p]$, i.e., there exists $\eta' \in
\Z[\delta_p]$ such that $$\eta\eta'=\eta'\eta=p \in
\Z[\delta_p].$$

 By a theorem of Ribet \cite{Ribet2} the
$\Z_p$-Tate module $T_p(J^{(f,p)})$ is a free
$\Z_p[\delta_p]$-module of rank $2g/(p-1)=n-1$ if $p$ does not
divide $n$ and $n-2$ if $p$ does. Let $J^{(f,p)}(\eta)$ be the
kernel of $\eta$ in $J^{(f,p)}(K_a)$. Clearly, $J^{(f,p)}(\eta)$ is
killed by multiplication by $p$, i.e.,  it may be viewed as a
$\F_p$-vector space. It follows from Ribet's theorem that
$$\dim_{\F_p}J^{(f,p)}(\eta)=\frac{2g}{p-1}.$$
In addition, $J^{(f,p)}(\eta)$ carries a
natural structure of Galois module. Notice that
$$\eta'(J^{(f,p)}_p) \subset J^{(f,p)}(\eta)$$
where $J^{(f,p)}_p$ is the kernel of multiplication by $p$ in $J^{(f,p)}(K_a)$.

Let $\Lambda$ be the centralizer of $\delta_p$ in $\End(J^{(f,p)})$. Clearly, $\Lambda$ commutes with $\eta$ and $\eta'$.
It is also clear that the subgroup $J^{(f,p)}(\eta)$ is $\Lambda$-stable. This observation leads to a natural homomorphism
$$\kappa:\Lambda \to \End_{\F_p}(J^{(f,p)}(\eta)).$$
I claim that its kernel coincides with $\eta\Lambda$.
Indeed, assume that $u(J^{(f,p)}(\eta))=\{0\}$ for some $u\in\Lambda$.
This implies that $u\eta'=\eta'u$ kills $J^{(f,p)}_p$. This implies, in turn,
 that there exists $v \in \End(J^{(f,p)})$ such that
$$u\eta'=\eta'u=pv=vp.$$
Clearly, $v$ commutes with $\eta$ and therefore with $\delta_p=1-\eta$, i.e., $v\in \Lambda$.
Since $p=\eta\eta'=\eta'\eta$,
$$u\eta'=v\eta\eta'$$
and therefore $u=v\eta$. On the other hand, it is clear that $\eta\Lambda=\Lambda\eta$ kills $J^{(f,p)}(\eta)$. Therefore the natural map
$$\Lambda/\eta\Lambda \to \End_{\F_p}(J^{(f,p)}(\eta))$$
is an embedding; further we will identify $\Lambda/\eta\Lambda$ with its image in $\End_{\F_p}(J^{(f,p)}(\eta))$.

\begin{thm}[Prop. 6.2 in \cite{Poonen}, Prop. 3.2 in \cite{SPoonen}]
\label{kereta}
There is a canonical isomorphism of the $\Gal(K)$-modules
$$ J^{(f,p)}(\eta) \cong V_{f,p}.$$
\end{thm}

\begin{rem}
\label{Galf}
Clearly, the natural homomorphism $\Gal(K) \to
\Aut_{\F_p}(V_{f,p})$ coincides with the composition
 $$\Gal(K)\twoheadrightarrow
\Gal(f)\subset \Perm(\R_f) \subset \Aut((\F_p^{\R_f})^{00})= \Aut_{\F_p}(V_{f,p}).$$
\end{rem}

The following assertion is an immediate corollary of Lemma \ref{Anp}.

\begin{lem}
\label{Anpf}
Assume that  $\Gal(f)= \Sn$ or $\An$.
If $n \ge 5$ then the $\Gal(f)$-module $V_{f,p}$ is absolutely simple.
\end{lem}

\begin{thm}
\label{notss}
Assume that $p>2$ and $n \ge 4$. Let $\Lambda_{\Q}=\Lambda\otimes\Q$
 be the centralizer of $\Q(\delta_p)$ in $\End^0(J^{(f,p)})$.
 Then $\Lambda_{\Q}$ could not be a central simple $\Q(\delta_p)$-algebra of dimension $(2g/(p-1))^2$ where $g$ is genus of $C_{f,p}$.
\end{thm}

\begin{proof}
Assume that $\Lambda_{\Q}$ is a central simple $\Q(\delta_p)$-algebra of dimension $(2g/(p-1))^2$.
We need to arrive to a contradiction. We start with the following statement.

\begin{lem}
\label{cm} Assume that $\Lambda_{\Q}$ is a central simple
$\Q(\delta_p)$-algebra of dimension $(2g/(p-1))^2$. Then there
exist a $(p-1)/2$-dimensional abelian variety $Z$ over $K_a$, a
positive integer $r$, an embedding $$\Q(\zeta_p)\cong \Q(\delta_p)
\hookrightarrow \End^0(Z)$$ and an isogeny $\phi: Z^r \to
J^{(f,p)}$ such that the induced isomorphism $$\Mat_r(\End^0(Z))
=\End^0(Z^r) \cong \End^0(J^{(f,p)}), \quad u \mapsto \phi
u\phi^{-1}$$ maps identically $$\Q(\delta_p) \subset \End^0(Z)
\subset \Mat_r(\End^0(Z)) =\End^0(Z^r)$$ onto $$\Q(\delta_p)
\subset \End^0(J^{(f,p)}).$$ (Here $\End^0(Z) \subset
\Mat_r(\End^0(Z))$ is the diagonal embedding.) In particular, $Z$
and $J^{(f,p)}$ are abelian varieties of CM-type over $K_a$.
\end{lem}

\begin{proof}[Proof of Lemma \ref{cm}]
Clearly, there exist a positive integer $r$  and a  central
division algebra $H$ over $\Q(\delta_p)\cong \Q(\zeta_p)$
 such that $\Lambda_{\Q}\cong \Mat_r(H)$. This imples that there exist an abelian variety $Z$ over $K_a$ with
$$\Q(\delta_p)\subset H \subset \End^0(Z)$$
 and  an isogeny
$\phi: Z^r \to J^{(f,p)}$
such that the induced isomorphism
$\End^0(Z^r) \cong \End^0(J^{(f,p)})$
maps identically
$$\Q(\delta_p)\subset\End^0(Z)\subset\End^0(Z^r)$$ onto
$\Q(\delta_p) \subset \End^0(J^{(f,p)})$.
We still have to check that
$$2\dim(Z)=p-1.$$
In order to do that let us put $g'=g/r$. Then $\dim_{\Q(\delta_p)}(H)= (\frac{2g'}{p-1})^2$ and therefore $\dim_{\Q}(H)= \frac{(2g')^2}{p-1}$. Since $H$ is a division algebra and $\fchar(K_a)=0$, the number $\frac{(2g')^2}{p-1}$ must divide $2\dim(Z)=2g'$. This means that $2g'$  divides $p-1$. On the other hand,
$$\Q(\delta_p) \subset H \subset \End^0(Z).$$
This implies that $p-1=[\Q(\delta_p):\Q]$ divides $2\dim(Z)$ and therefore $2\dim(Z)=p-1.$
\end{proof}

Now let us return to the proof of Theorem \ref{notss}.
Recall that $n \ge 4$. We write $\Omega^1(X)$ for the space of differentials of first kind for any smooth projective variety $X$ over $K_a$.
 Clearly, $\phi$ induces an isomorphism
$\phi^*: \Omega^1(J^{(f,p)}) \cong \Omega^1(Z^r)=\Omega^1(Z)^r$
which commutes with the natural actions of $\Q(\delta_p)$. Since
$\dim(Z)=\frac{p-1}{2}$, we have
$\dim_{K_a}(\Omega^1(Z))=\frac{p-1}{2}$. Therefore, the induced
$K_a$-linear automorphism $$\delta_p^*:\Omega^1(Z) \to
\Omega^1(Z)$$ has, at most, $\frac{p-1}{2}$ distinct eigenvalues.
Clearly, the same is true for the action of $\delta_p$ in
$\Omega^1(Z)^r$. Since $\phi$ commutes with $\delta_p$, the
induced $K_a$-linear automorphism $$\delta_p^*:
\Omega^1(J^{(f,p)})\to \Omega^1(J^{(f,p)})$$
 has, at most, $\frac{p-1}{2}$ distinct eigenvalues.

 On the other hand, let $P_0$ be one of the $\delta_p$-invariant points
(i.e., a ramification point for $\pi$) of $C_{f,p}(K_a)$. Then
$$\tau: C_{f,p} \to J^{(f,p)}, \quad P\mapsto \cl((P)-(P_0))$$ is
an embedding of $K_a$-algebraic varieties and it is well-known
that the induced map $$\tau^*: \Omega^1(J^{(f,p)}) \to
\Omega^1(C_{f,p})$$ is a $K_a$-linear isomorphism obviously
commuting with the actions of $\delta_p$. (Here $\cl$ stands for
the linear equivalence class.) This implies that $\delta_p$ has,
at most, $\frac{p-1}{2}$ distinct eigenvalues in
$\Omega^1(C_{f,p})$.

 One may easily
check that $\Omega^1(C_{f,p})$ contains differentials $dx/y^i$
for all positive integers $i<p$ satisfying $n i \ge (p+1)$ if $p$
does not divide $n$ (\cite{Koo}, Th. 3 on p. 403; see also
\cite{Towse}, Prop. 2 on p. 3359). Since $n \ge 4$ and $p \ge 3$,
we have $n i \ge (p+1)$ for all $i$ with
 $\frac{p-1}{2} \le i \le p-1$. Therefore  the differentials $dx/y^i \in \Omega^1(C_{f,p})$ for all $i$ with
 $\frac{p-1}{2} \le i \le p-1$;
 clearly, they all are eigenvectors of $\delta_p$ with eigenvalues $\zeta^{-i}$
 respectively.
(Recall that $\zeta\in K_a$ is a primitive $p$th root of unity
and $\delta_p$ is defined in \S 1 by $(x,y) \mapsto (x,\zeta y)$.)
Therefore $\delta_p$ has in $\Omega^1(C_{f,p})$, at least,
$\frac{p+1}{2}$ distinct eigenvalues. Contradiction.

Now assume that $p$ divides $n$. Then  $n \ge 5$. By Remark
\ref{nondiv}, $C_{f,p}$ is birationally isomorphic over $K_a$ to
a curve $C_1=C_{h_1,p}: y_1^p=h_1(x_1)$ where $h_1(x_1) \in
K_a[x_1]$ is a separable polynomial of degree $n-1$;
 in addition, one may choose this isomorphism in such a way that it commutes
with the actions of $\delta_p$ on $C_{f,p}$ and $C_{h_1,p}$. This
implies that $\delta_p$ has, at most, $\frac{p-1}{2}$ distinct
eigenvalues in $\Omega^1(C_{h_1,p})$.

On the other hand, $n-1 \ge 4$  and $n-1$ is not divisible by $p$.
Recall that $\deg(h_1)=n-1$. We  conclude, as above,
 that for all $i$ with $\frac{p-1}{2} \le i \le p-1$ the
 differentials $dx_1/y_1^i \in \Omega^1(C_{h_1,p})$.
Now, the same arguments as in the case of $p$ not dividing $n$
 lead to a contradiction.
\end{proof}

\begin{thm}
\label{handyp3}
Suppose $n \ge 4$ and $p>2$.
Assume that
$\Q(\delta_p)$ is a maximal commutative subalgebra
in $\End^0(J^{(f,p)})$.
Then:
\begin{enumerate}
\item[(i)]
The center $\CC$ of $\End^0(J^{(f,p)})$ is a CM-subfield of $\Q(\delta_p)$;
\item[(ii)]
If $p$ is a Fermat prime then $\End^0(J^{(f,p)})=\Q(\delta_p)
\cong\Q(\zeta_p)$ and therefore
 $\End(J^{(f,p)})=\Z[\delta_p] \cong \Z[\zeta_p]$.
\end{enumerate}
\end{thm}

\begin{proof}
Clearly,
${\CC} \subset \Q(\delta_p)$. Since $\Q(\delta_p)$ is a CM-field,
 $\CC$ is either a totally real field or a CM-field.
If $p$ is a Fermat prime then each subfield of $\Q(\delta_p)$
(distinct from $\Q(\delta_p)$ itself) is totally real. Therefore,
(ii) follows from (i).

In order to prove (i), let us assume that $\CC$ is totally real. We are going to arrive to a contradiction which proves (i).
Replacing, if necessary, $K$ by its subfield finitely generated over the
rationals, we may assume that $K$ (and therefore $K_a$) is isomorphic to
a subfield of the field $\C$ of complex numbers.
Since the center $\CC$ of  $\End^0(J^{f,p})$ is totally real,
the Hodge group of $J^{(f,p)}$ must be semisimple.
This implies that the pair $(J^{(f,p)},\Q(\delta_p))$ is of {\sl Weil type} (\cite{MZ}),
i.e., $\Q(\delta_p)$ acts on $\Omega^1(J^{(f,p)})$
in such a way that for each  embedding
$$\sigma: \Q(\delta_p) \hookrightarrow \C$$
the corresponding multiplicity
$$n_{\sigma}=
\frac{\dim(J^{(f,p)})}{[\Q(\delta_p):\Q]} .$$

Now assume that $p$ does not divide $n$. We have
$$\frac{\dim(J^{(f,p)})}{[\Q(\delta_p):\Q]}
=\frac{g(C_{f,p})}{p-1}=\frac{(n-1)}{2}$$ and therefore
$$n_{\sigma}=\frac{(n-1)}{2}.$$ Since the multiplicity
$n_{\sigma}$ is always an integer, $n$ is odd. Therefore $n \ge
5$. Let us consider the embedding $\sigma$ which sends $\delta_p$
to $\zeta$. Elementary calculations (\cite{Koo}, Th. 3 on p. 403)
show that for all integers $i$ with
 $$0 \le i \le n-1-\frac{(n+1)}{p}$$
the differentials $x^idx/y^{p-1} \in \Omega^1(C_{f,p})$;
 clearly, they constitute a set of $K_a$-linearly independent eigenvectors of $\delta_p$
 with eigenvalue $\zeta$. In light of the $\delta_p$-equivariant isomorphism
$$\Omega^1(J^{(f,p)}) \to \Omega^1(C_{f,p}),$$
we conclude that
$$\frac{(n-1)}{2}=n_{\sigma} \ge [n-1-\frac{(n+1)}{p}]+1.$$
This implies that $\frac{(n-1)}{2}>n-1-\frac{(n+1)}{p}$. It follows easily that
$n< \frac{p+2}{p-2} \le 5$ and therefore $n<5$.
This gives us the desired contradiction when $p$ does not divide $n$.

Now assume that $p$ divides $n$. Then $n \ge 5$ and $n-1 \ge 4$.
Again, as in the proof of Theorem \ref{notss}, the usage of Remark \ref{nondiv} allows us to apply the already proven case (when $p$ does not divide $n-1$)  to $C_{h_1,p}$ with $\deg(h_1)=n-1$.
\end{proof}

\begin{rem}
\label{hodge} Let us keep the notations and assumptions of
Theorem \ref{handyp3}. Assume, in addition that $p=3$. Then
$\Q(\delta_3)=\Q(\zeta_3)$ is an imaginary quadratic field and
there are exactly two embeddings $\Q(\delta_3) \hookrightarrow
K_a$ which, of course, are complex-conjugate. In this case one
could compute explicitly the corresponding multiplicities.

Indeed, first assume that $3$ does not divide $n$. Then
 $n=3k-e$ for some $k,e \in \Z$ with $3>e>0$.
Since $n \ge 4>3$, we have $k \ge 2$.
By Prop. 2 on p. 3359 of \cite{Towse}, the set
 $$\{x^i dx/y, 0\le i < k-1; x^j dx/y^2, 0 \le j< 2k-1-[\frac{2e}{3}]\}$$
is a basis of $\Omega^1(C_{(f,3)})$. It follows easily that it is
an eigenbasis with respect to the action of $\delta_3$. This
implies easily that $\Q(\delta_3)$ acts on
$\Omega^1(J^{(f,3)})=\Omega^1(C_{(f,3)})$ with multiplicities
$k-1$ and $2k-1-[\frac{2e}{3}]$.

Assume now that $n=3k$ is divisible by $3$. Then
as in the proof of Theorem \ref{notss}, the usage of Remark \ref{nondiv} allows us to reduce
 the calculation of multiplicities to the case of $C_{h_1,3}$ with $\deg(h_1)=n-1$.
More precisely, we have $n=3k$ and $n-1=3k-1$, i.e., $e=1$ and $k
\ge 2$. This implies that $\Q(\delta_3)$ acts on
$\Omega^1(J^{(f,3)})=\Omega^1(J^{(h_1,3)})$ with  multiplicities
$k-1$ and $2k-1$.

It follows that if $ n=3k$ or $n=3k-1$ then
$$\dim(J^{(f,3)})=3k-2$$
 and the imaginary quadratic field $\Q(\delta_3)$ acts on
$\Omega^1(J^{(f,3)})$ with {\sl mutually prime} multiplicities $k-1$ and $2k-1$. Since $\Q(\delta_3)$ coincides with $\End^0(J^{(f,3)})$,
a theorem of Ribet (\cite{Ribet3}, Th. 3 on p. 526) implies that the Hodge group of $J^{(f,3)}$
  is {\sl as large as possible}. More precisely, let $K' \subset K$ be a
subfield which admits an embedding into $\C$ and such that
$f(x) \subset K'[x]$ (such a subfield always exists).
Then one may consider $C_{f,3}$ as a complex smooth projective curve and
 $J^{(f,3)}$ as a complex abelian variety, whose endomorphism algebra coincides with
  $\Q(\delta_3)\cong \Q(\zeta_3)=\Q(\sqrt{-3})$.
Then the Hodge group of
$J^{(f,3)}$ coincides with the corresponding unitary group
 of $H_1(J^{(f,3)}(\C),\Q)$ over  $\Q(\zeta_3)$.
In particular, all Hodge classes on all self-products of
$J^{(f,3)}$ could be presented as linear combinations of exterior
products of divisor classes. As was pointed out in \cite{MZD},
pp. 572--573, the same arguments work also for Tate classes if
say, $K'$ is a number field and all endomorphisms of $J^{(f,3)}$
are defined over $K'$. (If $\sqrt{-3} \in K'$ then all
endomorphisms of $J^{(f,3)}$ are defined over $K'$ because
$\End(J^{(f,3)})=\Z[\delta_3]$ and $\delta_3$ is defined over
$K'(\zeta_3)=K'(\sqrt{-3})$.)
\end{rem}

\section{Representation theory}
\begin{defn}
Let $V$ be a vector space over a field $F$, let $G$ be a group and
$\rho: G \to \Aut_F(V)$ a linear representation of $G$ in $V$.
 We
say that the $G$-module $V$ is {\sl very simple} if it enjoys the
following property:

If $R \subset \End_F(V)$ be an $F$-subalgebra containing the
identity operator $\I$ such that

 $$\rho(\sigma) R \rho(\sigma)^{-1} \subset R \quad \forall
 \sigma \in G$$
 then either $R=F\cdot \I$ or $R=\End_F(V)$.
\end{defn}

\begin{rem}
\label{image}
\begin{enumerate}
\item[(i)]
 Clearly, the $G$-module $V$ is very simple if and only if the corresponding $\rho(G)$-module $V$ is very simple. It is known (\cite{ZarhinM}, Rem. 2.2(ii)) that a very simple module is absolutely simple.

\item[(ii)]
If $G'$ is a subgroup of $G$ and the $G'$-module $V$ is very simple then the $G$-module $V$ is also very simple.
\end{enumerate}
\end{rem}

\begin{thm}
\label{veryfactor}
Suppose a field $F$, a positive integer $N$ and a group $H$ enjoy the following properties:

\begin{itemize}
\item
 $F$ is either finite or algebraically closed;
 \item
 $H$ is perfect, i.e., $H=[H,H]$;
\item
Each homomorphism from $H$ to $\SS_N$ is trivial;
\item
Let $N=ab$ be a factorization of $N$ into a product of two
positive integers $a$ and $b$. Then either each homomorphism from
$H$ to $\PGL_a(F)$ is trivial or each homomorphism from $H$ to
$\PGL_b(F)$ is trivial.
\end{itemize}
Then each absolutely simple $H$-module of $F$-dimension $N$ is
very simple. In other words, in dimension $N$ the properties of
absolute simplicity and supersimplicity over $F$ are equivalent.
\end{thm}

\begin{proof}
We may assume that $N>1$. Let $V$ be an absolutely simple
$H$-module of $F$-dimension $N$.  Let $R \subset \End_F(V)$ be an
$F$-subalgebra containing the identity operator $\I$ and such that

 $$u R u^{-1} \subset R \quad \forall u \in H.$$

Clearly, $V$ is a faithful $R$-module and $$u R u^{-1} = R \quad
\forall u \in H.$$

{\bf Step 1}. By Lemma 7.4(i) of \cite{ZarhinM}, $V$ is a {\sl semisimple} $R$-module.

{\bf Step 2}. The $R$-module $V$ is {\sl isotypic}. Indeed, let us
split the semisimple $R$-module $V$ into the direct sum $$V =V_1
\oplus \cdots \oplus V_r$$ of its isotypic components.
 Dimension arguments imply that $r \le \dim(V) = N$.
 It follows easily from the arguments of the previous step that for each isotypic component $V_i$ its image
  $s V_i$ is an isotypic $R$-submodule for each $s \in H$ and therefore is contained in some $V_j$.
  Similarly, $s^{-1}V_j$ is an isotypic submodule obviously containing $V_i$. Since $V_i$ is the isotypic component,
  $s^{-1}V_j=V_i$ and therefore $sV_i=V_j$.
  This means that $s$ permutes the $V_i$; since $V$ is $H$-simple, $H$ permutes them transitively.
 This gives rise to the homomorphism $H \to \SS_r$ which must be  trivial,
 since $r \le N$ and therefore  $\SS_r$ is a subgroup of $\SS_{N}$.
This means that  $s V_i=V_i$ for all $s \in H$ and $V=V_i$ is
isotypic.

{\bf Step 3}. Since $V$ is isotypic, there exist a simple
$R$-module $W$ and a positive integer $d$ such that $V \cong W^d$.
We have
    $$d \cdot \dim(W)=\dim(V)=N.$$
Clearly, $\End_R(V)$ is isomorphic to the matrix algebra
$\Mat_d(\End_R(W))$  of size $d$ over $\End_R(W)$.

Let us put
    $$k=\End_R(W).$$
Since $W$ is simple, $k$ is a finite-dimensional division algebra over $F$. Since $F$ is either finite or algebraically closed, $k$ must be a field. In addition, $k=F$ if $F$ is algebraically closed and $k$ is finite if $F$ is finite. We have
    $$\End_R(V) \cong \Mat_d(k).$$
Clearly, $\End_R(V) \subset \End_{F}(V)$ is stable under the adjoint action of $H$. This induces a homomorphism     $$\alpha: H \to\Aut_F(\End_R(V))=\Aut_F(\Mat_d(k)).$$
 Since $k$ is the center of $\Mat_d(k)$, it is stable under the action of $H$, i.e., we get a homomorphism $H\to\Aut(k/F)$, which must be trivial, since $H$
is perfect and $\Aut(k/F)$ is abelian. This implies that the
center $k$ of $\End_R(V)$ commutes with $H$. Since $\End_H(V)=F$, we have $k=F$. This implies that
$\End_R(V) \cong \Mat_d(F)$ and
$$\alpha: H \to \Aut_F(\Mat_d(F))=\GL(d,F)/F^*=\PGL_d(F)$$
 is trivial if and only if
$\End_R(V) \subset \End_H(V)=F \cdot \I$.
 Since $\End_R(V) \cong \Mat_d(F)$, $\alpha$ is trivial if and only if $d=1$, i.e. $V$ is an absolutely simple $R$-module.

 It follows from the Jacobson density theorem that
 $R \cong \Mat_m(F)$ with $dm=N$. This implies that $\alpha$ is trivial if and only if
$R \cong \Mat_{N}(F)$, i.e., $R=\End_{\F}(V)$.

The adjoint action of $H$ on $R$ gives rise to a homomorphism
    $$\beta: H \to \Aut_F(\Mat_m(F))=\PGL_m(F).$$
 Clearly, $\beta$ is trivial if and only if $R$ commutes with $H$, i.e. $R=F \cdot\I$.

 It follows that we are done if either $\alpha$ or $\beta$ is trivial. Now one has only to recall that $N=dm$.
\end{proof}

\begin{cor}
\label{superp}
Let $p$ be a prime, $V$ a vector space over $\F_p$
of finite dimension $N$.
 Let $H \subset \Aut(V)$ be a non-abelian
simple group. Suppose that the $H$-module $V$ is absolutely simple and $H$ is not isomorphic to a subgroup of
 $ \SS_{N}$. Then the $H$-module $V$ is very simple if one of the following conditions holds:

\begin{enumerate}
\item[(i)]
$N$ is a prime;
\item[(ii)]
$N=8$ or twice a prime. In addition, $H$ is not
isomorphic to $PSL_2(\F_p)$ and either $H$ is not isomorphic to $\A_5$ or $p$ is not congruent to $\pm 1$ modulo $5$;
\item[(iii)]
$\#(H) \ge ((p^{[\sqrt{N}]}-1)^{[\sqrt{N}]})/(p-1)$
\item[(iv)]
$\#(H) \ge (p^{N}-1)/(p-1)$.
\end{enumerate}
\end{cor}

\begin{proof}
Let us split $N$ into a product $N=ab$ of two positive integers $a$ and $b$. In the case (i) either $a$ or $b$ is $1$ and the target of the corresponding projective linear group $\PGL_1(\F_p)=\{1\}$. In the case (ii) either one of the factors is $1$ and we are done or one of the factors is $2$ and it suffices to check that each homomorphism from
$H$ to $\PGL_2(\F_p)$ is trivial. Since $H$ is simple, each
non-trivial homomorphism $\gamma:H \to \PGL_2(\F_p)$ is an injection, whose image lies in $\PSL_2(\F_p)$. In other words, $\gamma(H)$ is
a subgroup of $\PSL_2(\F_p)$ isomorphic to $H$. Since $H$ is not isomorphic to $\PSL_2(\F_p)$, the subgroup $\gamma(H)$ is proper and simple non-abelian.
 It is known (\cite{Suzuki}, Th. 6.25 on p. 412 and Th. 6.26 on p. 414)
that each proper simple non-abelian subgroup of $\PSL_2(\F_p)$ is isomorphic to $\A_5$ and such a subgroup exists if and only if  $p$ is congruent to $\pm 1$ modulo $5$. This implies that such $\gamma$ does not exist and settles the case (ii).
In order to do the case (iii) notice that one of the factors say, $a$ does not exceed $[\sqrt{N}]$. This implies easily that the order of $\GL_a(\F_p)$ does not exceed
$((p^{[\sqrt{N}]}-1)^{[\sqrt{N}]})$ and therefore the order of $\PGL_a(\F_p)$ does not exceed
$((p^{[\sqrt{N}]}-1)^{[\sqrt{N}]})/(p-1)$. Hence, the order of $H$ is strictly greater than the order of $\PGL_a(\F_p)$ and therefore there are no injective homomorphisms from $H$ to $\PGL_a(\F_p)$. Since $H$ is simple, each homomorphism from $H$ is either trivial or injective. This settles the case (iii). The case (iv) follows readily from the case (iii).

\end{proof}

\begin{cor}
\label{Anvery3}
Suppose $n \ge 5$ is an integer,
 $B$ is an $n$-element set. Suppose $p=3$.
Then the $\Alt(B)$-module $(\F_3^{B})^{00}$ is very simple.
\end{cor}

\begin{proof}
By Lemma \ref{Anp}, $(\F_3^{B})^{00}$ is absolutely simple and
$N=\dim_{\F_3}((\F_3^{B})^{00})$ is either $n-1$ or $n-2$. The
group $\Alt(B)\cong \An$ is a simple non-abelian group, whose
order $n!/2$ is greater than the order of ${\mathbf S}_{n-1}$ and
the order of ${\mathbf S}_{n-2}$. Therefore  each homomorphism
from $\Alt(B)$ to ${\mathbf S}_N$ is trivial. On the other hand,
one may easily check that $$n!/2 > 3^{n-1}/2 > (3^N-1)/(3-1)$$ for
all $n \ge 5$. Now one has only to apply Corollary
\ref{superp}(iv) to $H=\Alt(B)$ and $p=3$.
\end{proof}

\begin{cor}
\label{verybign}
Suppose $p>3$ is a prime, $n \ge 8$ is a positive integer,  $B$ is an  $n$-element set.
Then the $\Alt(B)$-module $(\F_p^{B})^{00}$ is very simple.
\end{cor}

\begin{proof}
Recall that $N=\dim_{\F_p}((\F_p^{B})^{00}$ is either $n-1$ or $n-2$. In both cases
$$[\sqrt{N}]-1 < [n/3].$$
Clearly, $\Alt(B) \cong \An$ is perfect and every homomorphism from $\Alt(B)$ to ${\mathbf S}_N$ is trivial.

We are going to deduce the Corollary from Theorem \ref{veryfactor} applied to $F=\F_p$ and $H=\Alt(B)$.  In order to do that let us consider a factorization $N=ab$ of $N$ into a product of two positive integers $a$ and $b$. We may assume that $a>1, b>1$ and say, $a \le b$. Then
$$a-1 \le [\sqrt{N}]-1 <[n/3].$$
Let
$$\alpha: \An \cong \Alt(B) \to \PGL_a(\F_p)$$
be a group homomorphism. We need to prove that $\alpha$ is trivial. Let $\bar{\F}_p$ be an algebraic closure  of $\F_p$. Since $\PGL_a(\F_p)\subset \PGL_a(\bar{\F}_p)$, it suffices to check that the composition
$$ \An \cong \Alt(B) \to \PGL_a(\F_p)\subset \PGL_a(\bar{\F}_p)$$
which we continue denote by $\alpha$, is trivial.

Let $$\pi:\tilde{\A}_n\twoheadrightarrow \An$$
 be the universal central extension of the perfect group $\An$. It is well-known that $\tilde{\A}_n$ is perfect and the kernel (Schur's multiplier) of $\pi$ is a cyclic group of order $2$, since $n \ge 8$. One could lift $\alpha$ to the homomorphism
$$\alpha': \tilde{\A}_n \to \GL_a(\bar{\F}_p).$$ Clearly,
$\alpha$ is trivial if and only if $\alpha'$ is trivial. In order
to prove the triviality of $\alpha'$, let us put $m=[n/3]$ and
notice that $\An$ contains a subgroup $D$ isomorphic to
$(\Z/3\Z)^m$ (generated by disjoint $3$-cycles). Let $D'$ be a
Sylow $3$-subgroup in $\pi^{-1}(D)$. Clearly, $\pi$ maps $D'$
isomorphically onto $D$. Therefore, $D'$ is a subgroup of
$\tilde{\A}_n$ isomorphic to $(\Z/3\Z)^m$.

Now, let us discuss the image and the kernel of $\alpha'$. First,
since $\tilde{\A}_n$ is perfect, its image lies in
$\SL_a(\bar{\F}_p)$, i.e., one may view $\alpha'$ as a
homomorphism from $\tilde{\A}_n$ to $\SL_a(\bar{\F}_p)$. Second,
the only proper normal subgroup in $\tilde{\A}_n$  is the kernel
of $\pi$. This implies that if $\alpha'$ is nontrivial then its
kernel meets $D'$ only at the identity element and therefore
$\SL_a(\bar{\F}_p)$ contains the subgroup $\alpha'(D')$
isomorphic to $(\Z/3\Z)^m$. Since $p \ne 3$,  the group
$\alpha'(D')$  is conjugate to an elementary $3$-group of
diagonal matrices in $\SL_a(\bar{\F}_p)$. This implies that $$m
\le a-1.$$ Since $m=[n/3]$, we get a contradiction which implies
that our assumption of the nontriviality of $\alpha'$ was wrong.
Hence $\alpha'$ is trivial and therefore $\alpha$ is also trivial.
\end{proof}

\begin{thm}
\label{Altvery}
Suppose $n \ge 5$ is a positive integer, $B$ is an $n$-element set, $p$ is a prime. Then the $\Alt(B)$module $(\F_p^{B})^{00}$ is very simple.
\end{thm}

\begin{proof}
The case of $p=2$ was proven in \cite{ZarhinM}, Ex. 7.2. The case of
$p=3$ was done in Corollary \ref{Anvery3}. So, we may assume that $p\ge 5$.
In light of Corollary \ref{verybign} we may assume that $n<8$, i.e.,
$5 \le n \le 7$.

Assume that $n \ne p$. Then $p$ does not divide $n$ and $n-1$
is either a prime or twice a prime. Therefore
$$N=\dim_{\F_p}((\F_p^{B})^{00})=n-1$$
is either a prime or twice a prime. Now the very simplicity of $(\F_p^{B})^{00}$ follows from   the cases (i) and (ii) of Corollary \ref{superp}.

Assume now that $n=p$. Then either $n=p=5$ or $n=p=7$. In both
cases $$N=\dim_{\F_p}((\F_p^{B})^{00})=n-2$$ is a prime. Now the
very simplicity of $(\F_p^{B})^{00}$ follows from  Corollary
\ref{superp}(i).

\end{proof}

\section{Jacobians and endomorphisms}
Recall that $K$ is a field of characteristic zero, $f(x)\in K[x]$ is a polynomial of degree $n \ge 5$ without multiple roots, $\R_f \subset K_a$ the set of its roots, $K(\R_f)$ its
splitting field, $$\Gal(f)=\Gal(K(\R_f)/K)\subset\Perm(\R_f).$$

\begin{rem}
\label{root} Assume that $\Gal(f)=\Perm(\R_f)$ or $\Alt(\R_f)$.
Taking into account that $\Alt(\R_f)$ is non-abelian simple,
$\Perm(\R_f)/\Alt(\R_f) \cong \Z/2\Z$ and $K(\zeta)/K$ is
abelian, we conclude that the Galois group of $f$ over $K(\zeta)$
is also either $\Perm(\R_f)$ or $\Alt(\R_f)$. In particular, $f$
remains irreducible over $K(\zeta)$. So, in the course of the
proof of main results from Introduction we may assume that $\zeta
\in K$.
\end{rem}

\begin{thm}
\label{handysup} Let $p$ be an odd prime and $\zeta \in K$. If
the $\Gal(f)$-module $(\F_p^{\R_f})^{00}$  is very simple then
$\Q(\delta_p)$ coincides with its own centralizer in
$\End^0(J^{(f,p)})$ and the center of $\End^0(J^{(f,p)})$ is a
CM-subfield of $\Q(\delta_p)$. In particular, if $p$ is a Fermat
prime then $\End^0(J^{(f,p)})=\Q(\delta_p)$ and
$\End(J^{(f,p)})=\Z[\delta_p]$.
\end{thm}

Combining Theorems  \ref{handysup}, Remark \ref{root}, Theorem \ref{Altvery} and Remark \ref{image}(ii),
we obtain the following statement.

\begin{cor}
\label{handyV}
Let $p$ be an odd prime.
If $f(x)\in K[x]$ is an irreducible polynomial of degree $n\ge 5$ and $\Gal(f)=\Sn$ or $\An$ then
$\Q(\delta_p)$ is a maximal commutative subalgebra in $\End^0(J^{(f,p)})$ and the center of $\End^0(J^{(f,p)})$ is a CM-subfield of $\Q(\delta_p)$. In particular, if $p$ is a Fermat prime then
$\End^0(J^{(f,p)})=\Q(\delta_p)$
and
$\End(J^{(f,p)})=\Z[\delta_p]$.
\end{cor}

\begin{proof}[Proof of Theorem \ref{handysup}]
Recall that $J^{(f,p)}$ is a $g$-dimensional
abelian variety defined over $K$.

Since $J^{(f,p)}$ is defined over $K$, one may associate with every $u\in \End(J^{(f,p)})$ and $ \sigma \in \Gal(K)$ an endomorphism $^{\sigma}u \in \End(J^{(f,p)})$ such that $$^{\sigma}u(x)=\sigma
u(\sigma^{-1}x) \quad \forall x \in J^{(f,p)}(K_a).$$

Let us consider the centralizer $\Lambda$ of $\delta_p$ in
$\End(J^{(f,p)})$. Since $\delta_p$ is defined over $K$, we have
$^{\sigma}u \in \Lambda$ for all $u \in \Lambda$. Clearly,
$\Z[\delta_p]$ sits in the center of $\Lambda$ and the natural
homomorphism $$\Lambda\otimes \Z_p \to
\End_{\Z_p[\delta_p]}T_p(J^{(f,p)})$$ is an embedding. Here $
T_p(J^{(f,p)})$ is the $\Z_p$-Tate module of $J^{(f,p)}$ which is
a free $\Z_p[\delta_p]$-module of rank $\frac{2g}{p-1}$. Notice
that $$J^{(f,p)}(\eta)= T_p(J^{(f,p)})/\eta T_p(J^{(f,p)}).$$
Recall also that (Theorem \ref{kereta} and Remark \ref{Galf})
$$J^{(f,p)}(\eta)= (\F_p^{\R_f})^{00}$$ and $\Gal(K)$ acts on
$(\F_p^{\R_f})^{00}$  through $$\Gal(K)\twoheadrightarrow
\Gal(f)\subset \Perm(\R_f)\subset \Aut((\F_p^{\R_f})^{00}).$$
Since the $\Gal(f)$-module $(\F_p^{\R_f})^{00}$  is very simple,
the $\Gal(K)$-module $J^{(f,p)}(\eta)$ is also very simple,
thanks to Remark \ref{image}(i). On the other hand, if an
endomorphism $u \in \Lambda$ kills
$J^{(f,p)}(\eta)=\ker(1-\delta_p)$ then one may easily check that
there exists a unique $v \in \End(J^{(f,p)})$ such that $u= v
\cdot \eta$. In addition, $v \in \Lambda$. This implies that the
natural map $$\Lambda \otimes_{\Z[\delta_p]}\Z[\delta_p]/(\eta)
\to \End_{\F_p}(J^{(f,p)}(\eta))$$ is an embedding. Let us denote
by $R$ the image of this embedding. We have $$R:=\Lambda/\eta
\Lambda= \Lambda\otimes\Z_p/\eta \Lambda \otimes\Z_p \subset
\End_{\F_p}(J^{(f,p)}(\eta)).$$ Clearly, $R$ contains the
identity endomorphism and is stable under the conjugation via
Galois automorphisms. Since the $\Gal(K)$-module
$J^{(f,p)}(\eta)$ is very simple,
 either $R=\F_p \cdot \I$ or
 $R=\End_{\F_p}(J^{(f,p)}(\eta))$.
If $\Lambda/\eta\Lambda =R=\F_p \cdot \I$
 then $\Lambda$ coincides with $\Z[\delta_p]$.
This means that $\Z[\delta_p]$ coincides with its own centralizer
in $\End(J^{(f,p)})$ and therefore
$\Q(\delta_p)$ is a maximal commutative subalgebra in $\End^0(J^{(f,p)})$.

If $\Lambda/\eta\Lambda=R=\End_{\F_p}(J^{(f,p)}(\eta))$ then, by Nakayama's Lemma,
$$\Lambda\otimes \Z_p =
\End_{\Z_p[\delta_p]}T_p(J^{(f,p)})\cong
\Mat_{\frac{2g}{p-1}}(\Z_p[\delta_p]).$$
This implies easily that
the $\Q(\delta_p)$-algebra
$\Lambda_{\Q}=\Lambda\otimes\Q \subset \End^0(X)$ has
dimension $(\frac{2g}{p-1})^2$ and its center
has dimension $1$. This means that $\Lambda_{\Q}$ is a central
 $\Q(\delta_p)$-algebra of dimension
$(\frac{2g}{p-1})^2$.
 Clearly, $\Lambda_{\Q}$ coincides with the centralizer of
$\Q(\delta_p)$ in $\End^0(J^{(f,p)})$. Since $\delta_p$ respects the theta
 divisor on the jacobian $J^{(f,p)}$, the algebra $\Lambda_{\Q}$ is stable
under the corresponding Rosati involution and therefore is
semisimple as a $\Q$-algebra. Since its center is the field
$\Q(\delta_p)$, the $\Q(\delta_p)$-algebra $\Lambda_{\Q}$ is
central simple and has dimension $(\frac{2g}{p-1})^2$. By Theorem
\ref{notss}, this cannot happen. Therefore $\Q(\delta_p)$ is a
maximal commutative subalgebra in $\End^0(J^{(f,p)})$.
\end{proof}

{\bf Proof of main results}.
Clearly, Theorem \ref{endo} follows readily from Corollary \ref{handyV}.  Theorem \ref{trigonal} follows readily from Corollary \ref{handyV}   combined with Remark \ref{hodge}.


\begin{thebibliography}{99}
\bibitem{Atlas} J. H. Conway, R. T. Curtis, S. P. Norton, R. A. Parker, R. A. Wilson, Atlas of finite groups. Clarendon Press, Oxford, 1985.

\bibitem{CR} Ch. W. Curtis, I. Reiner, Representation theory of finite groups
and associative algebras. Interscience Publishers, New York London 1962.
\bibitem{Glasgow} H. K. Farahat, {\sl On the natural representation of the symmetric group}. Proc. Glasgow Math. Association {\bf 5} (1961-62), 121--136.


\bibitem{Isaacs} I. M.  Isaacs, Character theory of finite groups.
Academic Press, New York San Francisco London, 1976.


\bibitem{AtlasB} Ch. Jansen, K. Lux, R. Parker, R. Wilson,  An Atlas
of Brauer characters. Clarendon Press, Oxford, 1995.


\bibitem{Koo} J. K. Koo, {\em On holomorphic differentials of some algebraic function field of one variable over} C. Bull. Austral. Math. Soc. {\bf 43} (1991), 399--405.


\bibitem{MZD} B. Moonen, Yu. G. Zarhin, {\em Hodge and Tate classes on simple abelian fourfolds}. Duke Math. J. {\bf 77} (1995), 553--581.

\bibitem{MZ} B. Moonen, Yu. G. Zarhin, {\em Weil classes on abelian varieties}. J. reine angew. Math. {\bf 496} (1998), 83--92.

\bibitem{MumfordAV} D. Mumford, Abelian varieties, Second edition. Oxford University Press, London, 1974.

\bibitem{Ribet2} K. Ribet, {\em Galois action on division points of Abelian varieties with real
multiplications}. Amer. J. Math. 98 (1976), 751--804.
\bibitem{Ribet3} K. Ribet, {\em Hodge classes on certain abelian varieties}. Amer. J. Math. {\bf 105} (1983), 523--538.

\bibitem{Passman} D. Passman, Permutation groups. W. A. Benjamin, Inc., New York-Amsterdam, 1968.

\bibitem{Poonen} B. Poonen, E. Schaefer, {\em Explicit descent for Jacobians of cyclic covers of the projective line}. J. reine angew. Math. {\bf 488} (1997), 141--188.

\bibitem{SPoonen} E. Schaefer, {\em Computing a Selmer group of a Jacobian using functions on the curve}. Math. Ann. {\bf 310} (1998), 447--471.

\bibitem{Serre} J.-P. Serre, Topics in Galois Theory. Jones and Bartlett Publishers, Boston-London, 1992. 163--176;

\bibitem{SerreRep} J.-P. Serre, Linear representations of finite groups. Springer-Verlag, 1977.



\bibitem{Suzuki} M. Suzuki, Group Theory I. Springer-Verlag, 1982.

\bibitem{Towse} C. Towse, {\em Weierstrass points on cyclic covers of the projective line}. Trans. AMS {\bf 348} (1996), 3355-3377.



\bibitem{Wagner} A. Wagner, {\em The faithful linear representations of least degree of} $\Sn$ {\em and} $\An$ {\em over a field of odd characteristic}.
 Math. Z. {\bf 154} (1977), 103--114.

\bibitem{Zarhin} Yu. G. Zarhin, {\em Hyperelliptic jacobians without complex multiplication}. Math. Res. Letters {\bf 6}(2000), 123--132.

\bibitem{ZarhinM} Yu. G. Zarhin, {\em Hyperelliptic jacobians and modular representations}, http://xxx.lanl.gov/abs/math.AG/0003002, to appear in  ``Moduli of abelian varieties" (C. Faber, G. van der Geer, F. Oort, eds.), Birkh\"auser.
\end{thebibliography}
\end{document}